\documentclass[twocolumn,psfig]{article}
\makeatletter
\def\@normalsize{\@setsize\normalsize{12pt}\xpt\@xpt
\abovedisplayskip 10pt plus2pt minus5pt\belowdisplayskip \abovedisplayskip
\abovedisplayshortskip \z@ plus3pt\belowdisplayshortskip 6pt plus3pt
minus3pt\let\@listi\@listI}

\def\subsize{\@setsize\subsize{12pt}\xipt\@xipt}
\def\section{\@startsection {section}{1}{\z@}{24pt plus 2pt minus 2pt}
{12pt plus 2pt minus 2pt}{\large\bf}}
\def\subsection{\@startsection {subsection}{2}{\z@}{12pt plus 2pt minus 2pt}
{12pt plus 2pt minus 2pt}{\subsize\bf}}
\makeatother

\usepackage{amssymb, amsmath, amsthm, mathrsfs}
\usepackage{verbatim}
\usepackage[margin=1in]{geometry}

\ifpdf
  \usepackage[pdftex]{graphicx}
\else
  \usepackage{graphicx}
\fi

\ifpdf
  \DeclareGraphicsExtensions{.pdf, .jpg}
\else
  \DelcareGraphicsExtensions
  {.eps, .jpg}
\fi

\def\qed{\hfill $\Box$\vspace{0.3cm}}
\def\pf{\noindent{\bf Proof. }}
\newtheorem{lemma}{Lemma}[section]

\newtheorem{theorem}{Theorem}
\newtheorem{corollary}[lemma]{Corollary}

\begin{document}

\date{}

\title{\Large\bf On the Decycling Number of Bubble-sort Star Graphs
\footnote{This research is
supported by the Ministry of Science and Technology of Taiwan under the Grants
MOST104--2221--E--141--002--MY3.}}

\author{Yu--Zhe Liu$^1$, Shyue--Ming Tang$^2$ and Jou--Ming Chang$^3$
\\
{\small $^{1,3}$Institute of Information and Decision Sciences,}\\
{\small National Taipei University of Business, Taipei, Taiwan.}\\
{\small $^{2}$Department of Psychology and Social Work,}\\
{\small National Defense University, Taipei, Taiwan.}\\
{\small \{kevin155096,tang1119\}@gmail.com, spade@ntub.edu.tw} }

\maketitle

\thispagestyle{empty}
\subsection*{\centering Abstract}

{\em
Bubble-sort star graphs are a combination of star graphs and bubble sort graphs. They are bipartite graphs and also form a family of Cayley graphs.
The decycling number of a graph is the minimum number of vertices whose removal from the graph results in an acyclic subgraph. In this paper, we prove the decycling number $D(n)$ of an $n$-dimensional bubble-sort star graph for $n\leqslant 5$. We also show $D(n)$ satisfies the inequalities for $n\geqslant 6$:
\[ \frac{n!(2n-5)+2}{4n-8} \leqslant D(n) < \]
\begin{equation*}
\left \{
\begin{array}{ll}
\frac{n!}{2}-(\frac{n}{2})! & if\ n\ is\ even, \\
\frac{n!}{2}-2(\frac{n-1}{2})!+1 & otherwise.
\end{array}
\right.
\end{equation*}
}

\section{Introduction}

Let $G$ = ($V$,$E$) be a simple graph, with vertex set $V$ and edge set $E$. A subset $D \subset V(G)$ is called a {\em decycling set} (also called {\em feedback vertex set}) if the subgraph $G \setminus D$ is acyclic. The minimum cardinality of $D$ is called the {\em decycling number}
(or {\em feedback number})\cite{bau2002,beineke1997,erdos1965}. Finding the decycling number of a graph $G$ is equivalent to finding the greatest order of an induced forest of $G$ \cite{erdos1986}. The decycling set problem is very important in designing operating systems and VLSI chips. In 1972, Karp has proved that the decycling set problem is NP-complete \cite{karp1972}. The problem remains NP-complete on graphs of maximum degree four \cite{ueno1988}.

Determining the decycling number is difficult, but determining the bound of decycling number for some specific graphs is feasible. Focardi et al. gave bounds for the decycling number of hypercubes \cite{focardi2000}. Many results on the upper bound of the decycling number are proposed, like star graphs \cite{wang2003}, ($n$, $k$)-star graphs \cite{wang2012}, bubble sort graphs \cite{wang2015}, and so forth \cite{chang2004,gao2015,kuo2009,lien2014,lien2015,wang2005,xu2011}.

In this paper, we would like to establish an upper bound for the decycling number on a special class of interconnection networks, called {\em bubble-sort star graphs} (BS graphs for short), which was first proposed by Chou et al. \cite{chou1996} in 1996. BS graphs are constructed like a combination of star graphs and bubble sort graphs \cite{akers1987,akers1989,lakshmivarahan1993}. All of the three graph classes are bipartite graphs and belong to the family of Cayley graphs. A BS graph $BS_n$ ($n \geqslant 3$) has $n!$ vertices labeled with distinct permutations of $\{1,2,\ldots ,n\}$. The edge set of $BS_n$ is determined by a {\em transition set} (or called {\em generator}) $\mathcal{T}=\{(1,2),(1,3),\ldots, (1,n),(2,3),(3,4),\ldots, (n-1,n)\}$. For every $(i,j)\in \mathcal{T}$ ($i < j$), two vertices $u=(u_1u_2\ldots u_n)$ and $v=(v_1v_2\ldots v_n)$ are connected if $(u_i=v_j)$, $(u_j=v_i)$, and $(u_k=v_k)$ for $k \neq i,j$. The degree of every vertex in a BS graph is equivalent to the cardinality of the transition set. As a result, $BS_n$ is a ($2n-3$)-regular. Figure 1 shows the examples of $BS_3$ and $BS_4$.  In $BS_4$, vertex (1234) has five neighbors, i.e., (2134), (3214), (4231), (1324) and (1243). Some researchers think $\mathcal{T}$ as the edge set of a {\em transitive graph} that makes the label transition intuitive.

\begin{figure}[htb]
\begin{center}
\includegraphics[width=3.0 in]{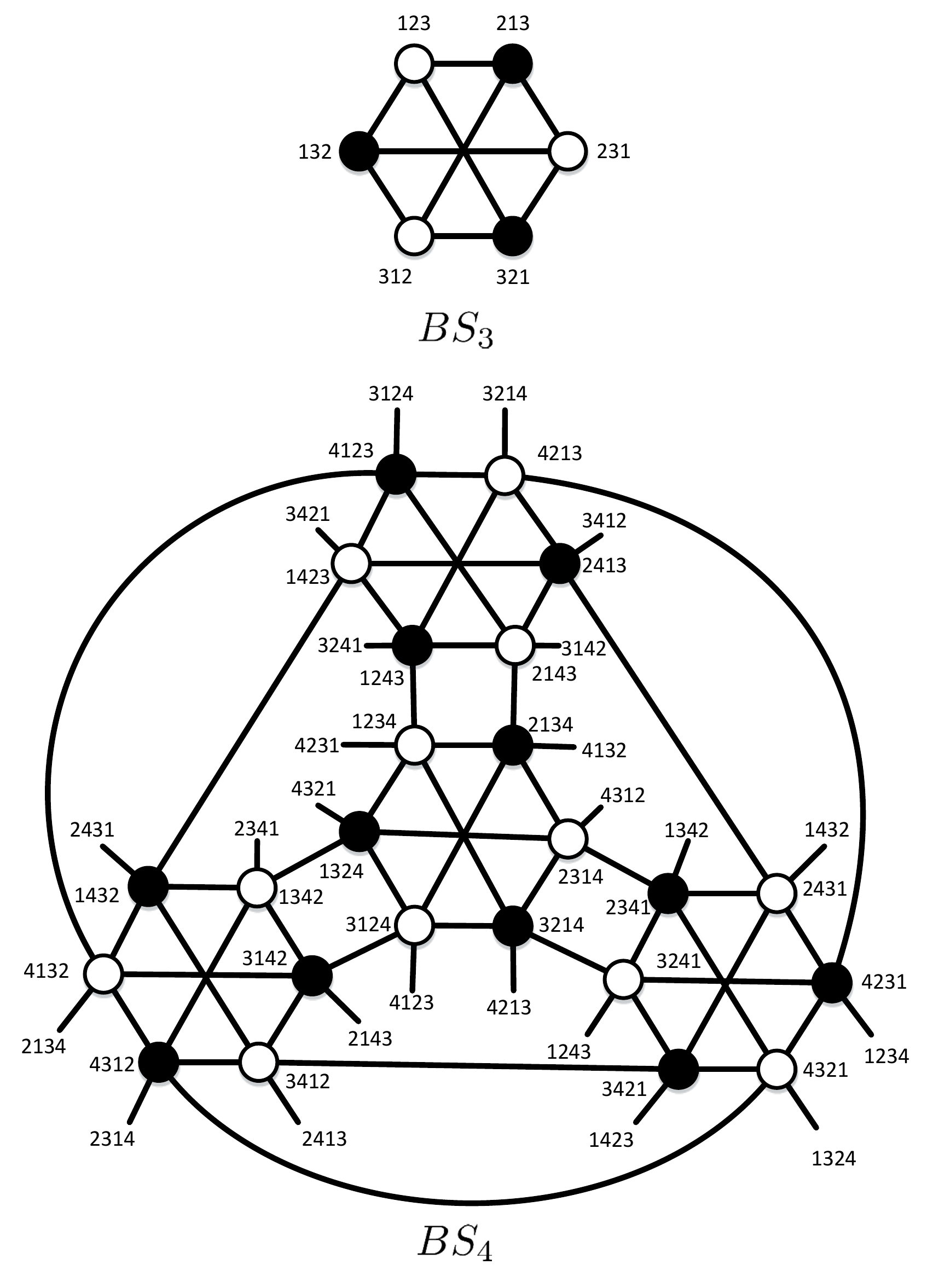} \\
\caption{The bubble-sort star graph samples, $BS_3$ and $BS_4$.}
\end{center}
\end{figure}

The BS graphs are intensively studied recently. In \cite{cai2015}, the authors showed that $BS_n$ is ($2n-5$)-fault-tolerant maximally local connected and is also ($2n-6$)-fault-tolerant one-to-many maximally local-connected. In \cite{guo2016a}, Guo et al. showed that for $BS_n$ the conditional diagnosability is $6n-15$ ($n > 5$) under the MM model and is $8n-21$ for $(n > 4)$ under the PMC model. In \cite{guo2016b}, they proved that $BS_n$ is $g$-connected for $1 \leqslant g \leqslant 3$. In \cite{wang2016}, Wang et al. proved that the 2-extra connectivity of $BS_n$ is $6n-15$ ($n > 4$) and the 2-extra connectivity of $BS_4$ is 8, the 2-extra diagnosability of $BS_n$ is $6n-13$ for $n > 4$ under the PMC model, and the 2-extra diagnosability of $BS_n$ is $6n - 13$ under the MM$^*$ model for $n > 5$. The pessimistic diagnosability of BS graphs is studied in \cite{gu2016}. In \cite{wang2017}, the authors prove that 2-good-neighbor connectivity of $BS_n$ is $8n-22$ for $n\geqslant 5$ and the 2-good-neighbor connectivity of $BS_4$ is 8. Further, the 2-good-neighbor diagnosability of $BS_n$ is $8n-19$ under the PMC model and MM$^*$ model for $n\geqslant 5$.

The remaining part of this paper is organized as follows. Section~2
presents the basic ideas of determining the decycling number of BS graphs. Section~3 gives the decycling number of $BS_n$ for $n\leqslant 5$. Section~4 gives a lower bound and an upper bound of the decycling number of $BS_n$ for $n\geqslant 6$. The last section contains our concluding remarks.

\section{Basic Ideas}

Let $(u_1u_2 \ldots u_n)$ be a permutation of the numerical set $\{1, 2, \ldots, n\}$. A numerical pair ($u_i, u_j$) is called an {\em inversion} if $i < j$ and $u_i > u_j$. In other words, an inversion in a permutation corresponds to a pair of numbers which are out of their natural order. For example, the permutation (41253) has four inversions, namely (4,1), (4,2), (4,3) and (5,3). In a BS graph, a vertex is called even (respectively, odd) if the number of inversions in its label is even (respectively, odd). In Figure 1, solid vertices are odd and blank vertices are even. The neighbors of an odd vertex must be even vertices, and vice versa. Since the BS graphs are bipartite and regular, the number of even and odd vertices is the same. Let $D(n)$ denote the decycling number of $BS_n$. Due to the bipartite property, there is an upper bound for $D(n)$. That is, $D(n)< \frac{n!}{2}$ since we can get an acyclic subgraph of $BS_n$ by removing all odd vertices (or even vertices, w. l. o. g.).

The {\em distance} $d(u,v)$ between two vertices $u$ and $v$ is defined as the length of a shortest path from $u$ to $v$.
An {\em independent set} (or called {\em stable set}) is a set of vertices in a graph $G$, no two of which are adjacent. A {\em distance-$k$ independent set} is a set $S$ of vertices in $G$ such that for any pair of vertices $u, v \in S$, the distance between $u$ and $v$ is at least $k$ in $G$. All even vertices in a BS graph form a independent set. A distance-4 independent set of odd vertices together with the even vertex set builds a acyclic subgraph of a BS graph. The basic ideas to reduce the upper bound of $D(n)$ is to remove a distance-4 independent set of odd vertices from the odd vertex set, such that the remained odd vertex set forms a smaller decycling set.

\section{The Decycling Numbers of $BS_n$ for $n\leqslant 5$}

In graph theory, The {\em diameter} of a connected graph $G$, denoted by $diam(G)$, is the greatest distance between any pair of vertices in the graph. Then we have $diam(BS_3)=2$, $diam(BS_4)=4$, and $diam(BS_5)=5$. Chou et al. give the following lemma.

\begin{lemma} {\rm \cite{chou1996}}
For $n>5$, $diam(BS_n)=\frac{3n}{2}-2$ if $n$ is even, $diam(BS_n)=\frac{3(n-1)}{2}$ otherwise.
\end{lemma}

A {\em rotation} on a permutation $(u_1,u_2,\ldots ,u_n)$ makes a new permutation $(v_1,v_2,\ldots, v_n)$ such that $v_i=u_{i+1}$ for $1\leqslant i\leqslant n-1$ and $v_n=u_1$. Let $R(u,k)$ denote the vertex whose label obtained from $u$ by performing $k$ rotations consecutively. Obviously, $R(u,n)=u$.

We show that the distance between vertices $u$ and $R(u,k)$ in $BS_n$ is always $n-1$ if $n$ is odd.

\begin{lemma}
For odd $n$, $d(u,R(u,k))=n-1$ in $BS_n$, where $k=1,2,\ldots,n-1$.
\end{lemma}

\pf Suppose $u=(u_1u_2\ldots u_n)$ be a vertex in $BS_n$. After $k$ rotations, $R(u,k)=(u_{k+1}u_{k+2}\ldots u_nu_1u_2\ldots u_k)$. A shortest path from $R(u,k)$ to $u$ is always taking the transition $(1,t)$, where $t$ is the first label number of the current vertex. Every transition gets exactly one number to its right target position, and no further transition needed for the number. More precisely, the $i$-th transition gets $u_t$ to its right label position, where $t=ki+1 ({\rm mod}\ n)$. After $n-1$ transitions, the path is completed. For example, vertex (45123) takes transition (1,4) to reach vertex (25143), then takes transition (1,2) to reach vertex (52143), and so forth. After four transitions, target vertex (12345) is arrived at.
\qed

Given $BS_3$ and $u=(123)$ as an example of Lemma 3.2, $R(u,1)=(231)$ and $R(u,2)=(312)$. Each vertex is distance 2 from each other .
According to the definition of distance-$k$ independent set, we have the following corollary.

\begin{corollary}
In $BS_n$ with odd $n$, vertices $u, R(u,1), R(u,2), \ldots, R(u,n-1)$ form a distance-$(n-1)$ independent set.
\end{corollary}

The followed lemma shows that the distance between vertices $u$ and $R(u,k)$ in $BS_n$ is greater than or equal to $n-1$ if $n$ is even.

\begin{lemma}
For even $n$, $d(u,R(u,k))=\frac{3n}{2}-2-|k-\frac{n}{2}|$ in $BS_n$, where $k=1,2,\ldots,n-1$.
\end{lemma}

\pf  The same transition policy for odd $n$ is applied to find a shortest path from $R(u,k)$ to $u$. There are $n-1$ necessary transitions. However, the right label position formula $t=ki+1 ({\rm mod}\ n)$ may get $t=1$ in case of even $n$. Every time the ``$t=1$'' case happens, one more transition $(1,s)$ is followed if there exists a label number $u_s$ that is not yet at its right position. The number $k$ of rotations determines how many times ``$t=1$'' case happens. A shortest path from $R(u,k)$ to $u$ gets ``$t=1$'' case $\frac{n}{2}-1-|k-\frac{n}{2}|$ times which is greater than or equal to 0. Consequently, the total distance from $R(u,k)$ to $u$ is $n-1+\frac{n}{2}-1-|k-\frac{n}{2}|=\frac{3n}{2}-2-|k-\frac{n}{2}|$.
\qed

In Lemma 3.4, the case $d(u,R(u,k))=n-1$ happens when $k=1$ and $k=n-1$. For example, $u=(123456)$ is a vertex in $BS_6$. Then, $R(u,1)=(234561)$ and $R(u,5)=(612345)$. We have $d(u,R(u,1))=d(u,R(u,5))=5$. As for other $R(u,k)$, we have $d(u,R(u,2))=d(u,R(u,4))=6$ and $d(u,R(u,3))=7$

The following corollary is derived from Lemma 3.1, Lemma 3.2 and Lemma 3.4.

\begin{corollary}
In $BS_n$ with $n>5$, $diam(BS_n)=d(u,R(u,\frac{n}{2}))$ if $n$ is even,
$diam(BS_n)=d(u,R(u,1))+\frac{n-1}{2}$ otherwise.
\end{corollary}

The following lemma gives a lower bound to the decycling number of a general graph.

\begin{lemma} {\rm \cite{beineke1997}}
If $D$ is a decycling set in a graph $G = (V, E)$ with maximum degree $\Delta$, the inequality must hold:
\[ |D| \geqslant \left \lceil  \frac{|E|-|V|+1}{\Delta -1} \right \rceil \ . \]
\end{lemma}

According to Lemma 3.6, we get the lower bound of $D(n)$ by giving $|E|=\frac{n!(2n-3)}{2}$, $|V|=n!$, and $\Delta=2n-3$.

\begin{corollary}
In $BS_n$, the lower bound of $D(n)$ is
\[D(n)\geqslant \frac{n!(2n-5)+2}{4n-8} \ .\]
\end{corollary}

The lower bounds of $D(3)$, $D(4)$ and $D(5)$ are 2, 10 and 51, respectively.
Since we have reduced the upper bounds to meet their lower bounds, the $D(n)$ is determined for $n\leqslant 5$.

\begin{theorem}
The decycling number of BS graphs, $D(3)=2$, $D(4)=10$, and $D(5)=51$.
\end{theorem}

\pf In $BS_3$, as shown in Figure 2, two gray vertices are a decycling set whose removal from $BS_3$ results in an acyclic subgraph. That is, $D(3)=2$. \\
In case of $BS_4$, only two out of the twelve odd vertices can be reserved. We choose any odd vertex $u=(u_1u_2u_3u_4)$, and then choose $R(u,2)=(u_3u_4u_1u_2)$. Since $d(u,R(u,2))=4$ and $R(u,2)$ must be an odd vertex, the lower bound of $D(4)$ is met. As shown in Figure 2, two solid vertices are reserved and ten gray vertices are removed. \\
In case of $BS_5$, we have to reserve nine odd vertices to meet the lower bound of $D(5)$. We can choose $\mathcal{O}_1=\{(u_1u_2u_3u_4u_5)$, $(u_2u_3u_4u_5u_1)$, $(u_3u_4u_5u_1u_2)$, $(u_4u_5u_1u_2u_3)$, $(u_5u_1u_2u_3u_4)\}$ as the first reserved odd vertex set, and $\mathcal{O}_2=\{(u_1u_2u_4u_5u_3)$, $(u_2u_3u_5u_1u_4)$, $(u_3u_4u_1u_2u_5)$, $(u_4u_5u_2u_3u_1)$, $(u_5u_1u_3u_4u_2)\}$ as the second reserved odd vertex set. By Corollary 3.3, set $\mathcal{O}_1$ is a distance-4 independent set. A shortest path between any two vertices of $\mathcal{O}_2$ can be constructed by using the same transition policy mentioned in Lemma 3.2. Thus set $\mathcal{O}_2$ is also a distance-4 independent set. In addition, every vertex in $\mathcal{O}_1$ exactly has two distance-2 neighbors and three distance-4 neighbors in $\mathcal{O}_2$, and vice versa. For example, vertex $(u_1u_2u_3u_4u_5)$ is two steps away from $(u_1u_2u_4u_5u_3)$ and $(u_5u_1u_3u_4u_2)$, and four steps away from $(u_2u_3u_5u_1u_4)$, $(u_3 u_4u_1u_2u_5)$ and $(u_4u_5u_2u_3u_1)$. Therefore, set $\mathcal{O}_1 \cup \mathcal{O}_2$ induces a 20-cycle with ten intermediate even vertices. As a result, if any odd vertex is removed from set $\mathcal{O}_1 \cup \mathcal{O}_2$, the cycle is broken and the lower bound of $D(5)$ is met.
\qed

\begin{figure}[htb]
\begin{center}
\includegraphics[width=3.0 in]{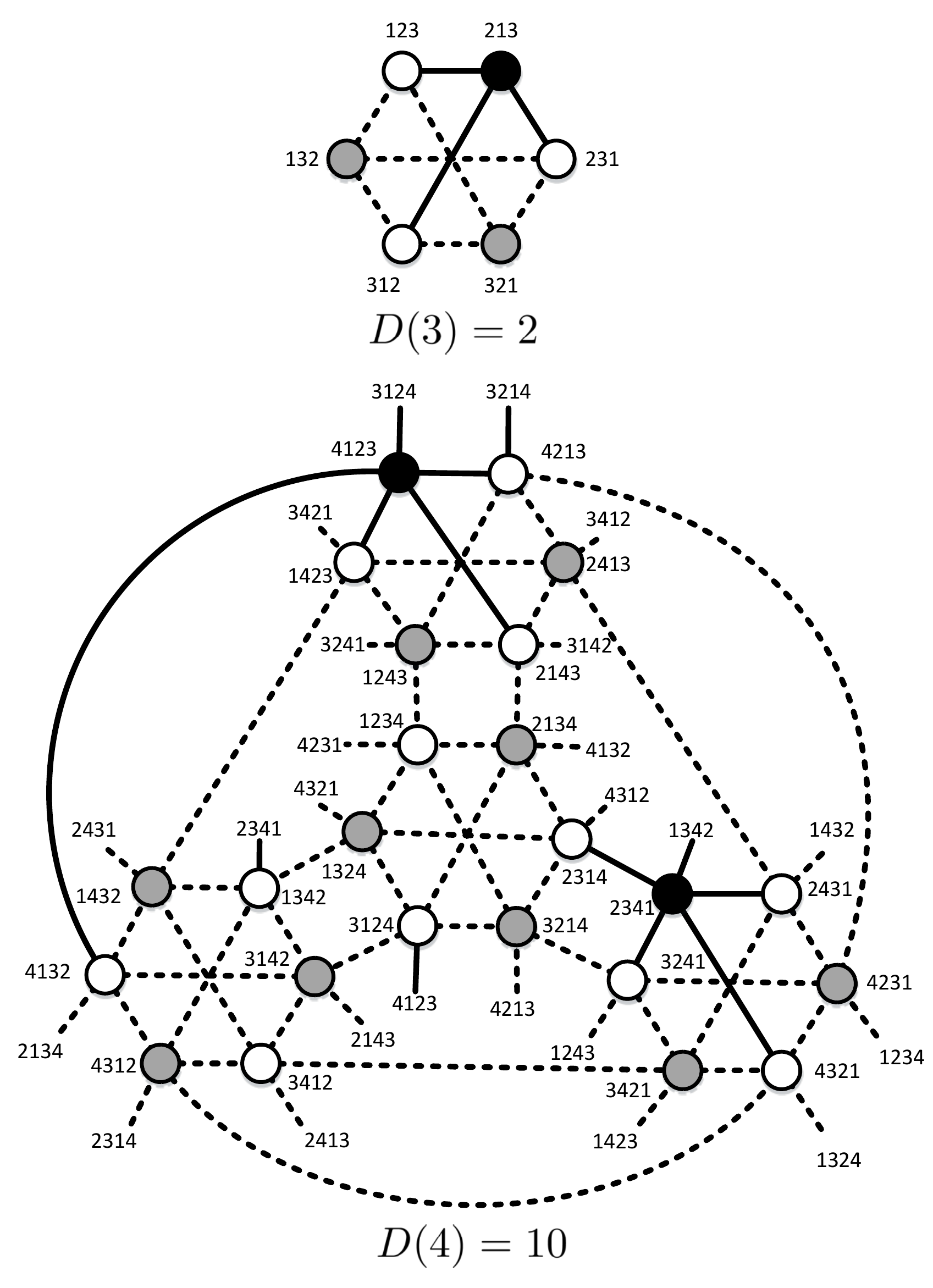} \\
\caption{The decycling number of $BS_3$ and $BS_4$.}
\end{center}
\end{figure}

\section{The Decycling Numbers of $BS_n$ for $n\geqslant 6$}

For $n\geqslant 6$, the upper bound of $D(n)$ is computed by finding as large as possible an odd vertex set that is distance-4 independent. The vertex set is reserved in a induced acyclic subgraph of $BS_n$ together with all even vertices to make the decycling set as small as possible.

Let $S=\{u_1, u_2, \ldots, u_n\}$ be a distinct number set and $n$ is even. The {\em consecutively paired number set} of $S$ is defined as $P(S)=\{u_1u_2, u_3u_3, \ldots, u_{n-1}u_n\}$.

Let $u=(u_1u_2\ldots u_n)$ be vertex in $BS_n$. If $n$ is even, then we have $A=\{u_1, u_2, \ldots, u_n\}$. The vertex set generated by permutations of $P(A)$ is called the {\em paired permutation vertex set} of $u$, denoted by $\mathcal{P}(u)$. For example, given a vertex $u=(123456)$ in $BS_6$, $\mathcal{P}(u)=\{$(123456), (125634), (341256), (345612), (561234), (563412)$\}$. The cardinality of $\mathcal{P}(u)$ is $\frac{n}{2}!$ if $n$ is even.

If $n$ is odd, then we have $A=\{u_1, u_2, \ldots, u_{n-1}\}$ and $B=\{u_1, u_3, u_4, \ldots, u_n\}$.
Note that vertices generated by permutations of $P(A)$ must end with $u_n$; while vertices generated by permutations of $P(B)$ must fix $u_2$ to the second position of the label. The paired permutation vertex set of $u$ includes vertices generated by permutations of $P(A)$ and permutations of $P(B)$. For example, given a vertex $u=(1234567)$ in $BS_7$, $\mathcal{P}(u)=\{$(1234567), (1256347), (3412567), (3456127), (5612347), (5634127), (1236745), (4251367), (4256713), (6271345), (6274513)$\}$. The cardinality of $\mathcal{P}(u)$ is $2\frac{n-1}{2}!-1$ if $n$ is odd.

The following lemma depicts the main result of this paper.

\begin{theorem}
An upper bounds of $D(n)$, for $n\geqslant 6$,  are
\[D(n) < \frac{n!}{2}-(\frac{n}{2})!\] if $n$ is even, or
\[D(n) < \frac{n!}{2}-2(\frac{n-1}{2})!+1\] if $n$ is odd.
\end{theorem}

\pf
If $u$ is an odd vertex, we prove the proposition holds. That is, vertex set $\mathcal{P}(u)$ is a distance-4 independent set and all vertices in the set are odd.\\
In case of $n$ is even, we prove this theorem by mathematical induction on $n$.
When $n=6$, for any odd vertex $u$ in $BS_6$, vertex set $\mathcal{P}(u)$ satisfies the proposition.
Suppose the proposition holds for $n=2k$. When $n=2k+2$, the vertex set $\mathcal{P}(u)$ can be partitioned into   $k+1$ vertex sets according to the position of new number pair ($u_{2k+1}u_{2k+2}$). The proposition holds in every vertex set, and for any two vertices from different vertex sets the distance between them is at least four.\\
In case of $n$ is odd, we prove this theorem by reducing $n$ to $n-1$. There are two induced subgraphs in $BS_n$ which have the same transition set isomorphic to that of $BS_{n-1}$. One subgraph is induced by vertices which have identical $u_n$ and another is induced by vertices which have identical $u_2$. Obviously, the vertex set $\mathcal{P}(u)$ obtained from the two subgraphs satisfy the proposition. Further, the distance between any two vertices from different vertex sets is at least four.\\
Consequently, an upper bound of $D(n)$ is \(\frac{n!}{2}-(\frac{n}{2})!\) if $n$ is even,
or \(\frac{n!}{2}-2(\frac{n-1}{2})!+1\) if $n$ is odd.
\qed

\section{Conclusion}

In this paper, we study the decycling numbers of BS graphs, and prove that $D(3)=2$, $D(4)=10$ and $D(5)=51$. As for $n \geqslant 6$, a lower bound of $D(n)$ is $\frac{n!(2n-5)+2}{4n-8}$ that is obtained directly from \cite{beineke1997}. An upper bound of $D(n)$ is $\frac{n!}{2}-(\frac{n}{2})!$ if $n$ is even,
or \(\frac{n!}{2}-2(\frac{n-1}{2})!+1\) if $n$ is odd. Further reduction of the upper bound of BS graphs is challenging. The study on the upper bound of decycling number in other Cayley graph classes is also challenging.

\end{document}